\documentclass[10pt,twoside]{article}
\usepackage{amsmath,amssymb}
\def\YEAR{\year}\newcount\VOL\VOL=\YEAR\advance\VOL by-1995
\def\firstpage{1}\def\lastpage{1000}
\def\received{}\def\revised{}
\def\communicated{}

\makeatletter
\def\magnification{\afterassignment\m@g\count@}
\def\m@g{\mag=\count@\hsize6.5truein\vsize8.9truein\dimen\footins8truein}
\makeatother

\font\eightrm=cmr8
\font\caps=cmcsc10                    % Theorem, Lemma etc
\font\Caps=cmcsc10 scaled \magstep1   % Title

\makeatletter
\@specialpagefalse\headheight=8.5pt
\def\DocMath{}
\renewcommand{\@evenhead}{%
    \ifnum\thepage>\lastpage\rlap{\thepage}\hfill%
    \else\rlap{\thepage}\slshape\leftmark\hfill{\caps\SAuthor}\hfill\fi}%
\renewcommand{\@oddhead}{%
    \ifnum\thepage=\firstpage{\DocMath\hfill\llap{\thepage}}%
    \else{\slshape\rightmark}\hfill{\caps\STitle}\hfill\llap{\thepage}\fi}%
\makeatother

\def\TSkip{\bigskip}
\newbox\TheTitle{\obeylines\gdef\GetTitle #1
\ShortTitle  #2
\SubTitle    #3
\Author      #4
\ShortAuthor #5
\EndTitle
{\setbox\TheTitle=\vbox{\baselineskip=20pt\let\par=\cr\obeylines%
\halign{\centerline{\Caps##}\cr\noalign{\medskip}\cr#1\cr}}%
	\copy\TheTitle\TSkip\TSkip%
\def\next{#2}\ifx\next\empty\gdef\STitle{#1}\else\gdef\STitle{#2}\fi%
\def\next{#3}\ifx\next\empty%
    \else\setbox\TheTitle=\vbox{\baselineskip=20pt\let\par=\cr\obeylines%
    \halign{\centerline{\caps##} #3\cr}}\copy\TheTitle\TSkip\TSkip\fi%
\centerline{\caps #4}\TSkip\TSkip%
\def\next{#5}\ifx\next\empty\gdef\SAuthor{#4}\else\gdef\SAuthor{#5}\fi%
\ifx\received\empty\relax
    \else\centerline{\eightrm Received: \received}\fi%
\ifx\revised\empty\TSkip%
    \else\centerline{\eightrm Revised: \revised}\TSkip\fi%
\ifx\communicated\empty\relax
    \else\centerline{\eightrm Communicated by \communicated}\fi\TSkip\TSkip%
\catcode'015=5}}\def\Title{\obeylines\GetTitle}
\def\Abstract{\begingroup\narrower
    \parskip=\medskipamount\parindent=0pt{\caps Abstract. }}
\def\EndAbstract{\par\endgroup\TSkip}

\long\def\MSC#1\EndMSC{\def\arg{#1}\ifx\arg\empty\relax\else
     {\par\narrower\noindent%
     2000 Mathematics Subject Classification: #1\par}\fi}

\long\def\KEY#1\EndKEY{\def\arg{#1}\ifx\arg\empty\relax\else
	{\par\narrower\noindent Keywords and Phrases: #1\par}\fi\TSkip}

\newbox\TheAdd\def\Addresses{\vfill\copy\TheAdd\vfill
    \ifodd\number\lastpage\vfill\eject\phantom{.}\vfill\eject\fi}
{\obeylines\gdef\GetAddress #1
\Address #2 
\Address #3
\Address #4
\EndAddress
{\def\xs{4.3truecm}\parindent=0pt
\setbox0=\vtop{{\obeylines\hsize=\xs#1\par}}\def\next{#2}
\ifx\next\empty % 1 address
     \setbox\TheAdd=\hbox to\hsize{\hfill\copy0\hfill}
\else\setbox1=\vtop{{\obeylines\hsize=\xs#2\par}}\def\next{#3}
\ifx\next\empty % 2 addresses
     \setbox\TheAdd=\hbox to\hsize{\hfill\copy0\hfill\copy1\hfill}
\else\setbox2=\vtop{{\obeylines\hsize=\xs#3\par}}\def\next{#4}
\ifx\next\empty\ % 3 addresses
     \setbox\TheAdd=\vtop{\hbox to\hsize{\hfill\copy0\hfill\copy1\hfill}
                \vskip20pt\hbox to\hsize{\hfill\copy2\hfill}}
\else\setbox3=\vtop{{\obeylines\hsize=\xs#4\par}}
     \setbox\TheAdd=\vtop{\hbox to\hsize{\hfill\copy0\hfill\copy1\hfill}
	        \vskip20pt\hbox to\hsize{\hfill\copy2\hfill\copy3\hfill}}
\fi\fi\fi\catcode'015=5}}\gdef\Address{\obeylines\GetAddress}

\begin{document}
%%%%% ------------- fill in your data below this line  -------------------
%%%%%    The following lines \Title ... \EndAddress must ALL be present
%%%%%    and in the given order.
\Title{Cohomology of arithmetic groups with infinite dimensional coefficient
spaces}
%%%%%    Put here the title. Line breaks will be recognized. 
\ShortTitle{COHOMOLOGY OF ARITHMETIC GROUPS...}
%%%%%    Running title for odd numbered pages, ONE line, please. 
%%%%%    If none is given, \Title will be used instead.          
\SubTitle 
%%%%%    A possible subtitle goes here.
\Author{Anton Deitmar \& Joachim Hilgert}
%%%%%    Put here name(s) of authors. Line breaks will be recognized.  
\ShortAuthor{Deitmar \& Hilgert}
%%%%%%   Running title for even numbered pages, ONE line, please. 
%%%%%%   If none is given, \Author will be used instead.          
\EndTitle
\Abstract 
{The cuspidal cohomology groups of arithmetic groups in certain infinite dimensional Modules are computed. As a result we get a simultaneous generalization of the Patterson-Conjecture and the Lewis-Correspondence.}
%%%%%    Put here the abstract of your manuscript.
%%%%%    Avoid macros and complicated TeX expressions, as this is
%%%%%    automatically translated and posted as an html file.
\EndAbstract
\MSC 11F75
%%%%%    2000 Mathematics Subject Classification: 
\EndMSC
\KEY 
%%%%%    Keywords and Phrases:     
\EndKEY
%%%%%    All 4 \Address lines below must be present. To center the last
%%%%%    entry, no empty lines must be between the following \Address
%%%%%    and \EndAddress lines.
\Address 
%%%%%    Address of first Author here
\Address
%%%%%    Address of second Author here etc.
\Address{Mathematisches Institut,  Auf der Morgenstelle 10, 72076 T\"ubingen,  Germany, \tt deitmar@uni-tuebingen.de}
\Address{Universit\"at Paderborn, Fakult\"at f\"ur Elektrotechnik, Informatik und
Mathematik, Warburger Str.  100, 33098 Paderborn, Germany.
\tt hilgert@math.uni-paderborn.de}
\EndAddress
%%
%%       Make sure the last tex command in your manuscript
%%       before the first \end{document} is the command  \Addresses
%%
%%---------------------Here the prologue ends---------------------------------
%%--------------------Here the manuscript starts------------------------------

\def \1{{\bf 1}}
\def \a{{{\mathfrak a}}}
\def \ad{{\rm ad}}
\def \al{\alpha}
\def \ar{{\alpha_r}}
\def \A{{\mathbb A}}
\def \Ad{{\rm Ad}}
\def \Aut{{\rm Aut}}
\def \b{{{\mathfrak b}}}
\def \bs{\backslash}
\def \B{{\cal B}}
\def \c{{\mathfrak c}}
\def \cent{{\rm cent}}
\def \C{{\mathbb C}}
\def \CA{{\cal A}}
\def \CB{{\cal B}}
\def \CC{{\cal C}}
\def \CD{{\cal D}}
\def \CE{{\cal E}}
\def \CF{{\cal F}}
\def \CG{{\cal G}}
\def \CH{{\cal H}}
\def \CHC{{\cal HC}}
\def \CL{{\cal L}}
\def \CM{{\cal M}}
\def \CN{{\cal N}}
\def \CP{{\cal P}}
\def \CQ{{\cal Q}}
\def \CO{{\cal O}}
\def \CS{{\cal S}}
\def \CT{{\cal T}}
\def \CU{{\cal U}}
\def\cusp{_{\rm cusp}}
\def \CV{{\cal V}}
\def \CW{{\cal W}}
\def \d{{\mathfrak d}}
\def \det{{\rm det}}
\def \df{\ \begin{array}{c} _{\rm def}\\ ^{\displaystyle =}\end{array}\ }
\def \diag{{\rm diag}}
\def \dist{{\rm dist}}
\def \ds{\displaystyle}
\def \End{{\rm End}}
\def \eps{\varepsilon}
\def \eqn{\begin{eqnarray*}}
\def \endeqn{\end{eqnarray*}}
\def \Ext{{\rm Ext}}
\def \F{{\mathbb F}}
\def \Fx{{\mathfrak x}}
\def \FX{{\mathfrak X}}
\def \g{{{\mathfrak g}}}
\def \ga{\gamma}
\def \Ga{\Gamma}
\def \Gal{{\rm Gal}}
\def \GL{{\rm GL}}
\def \h{{{\mathfrak h}}}
\def \H{{\mathbb H}}
\def \Hom{{\rm Hom}}
\def \im{{\rm im}}
\def \Im{{\rm Im}}
\def \Ind{{\rm Ind}}
\def \k{{{\mathfrak k}}}
\def \K{{\cal K}}
\def \l{{\mathfrak l}}
\def \la{\lambda}
\def \lap{\triangle}
\def \li{{\rm li}}
\def \La{\Lambda}
\def \m{{{\mathfrak m}}}
\def \mod{{\rm mod}}
\def \n{{{\mathfrak n}}}
\def \name{\bf}
\def \Mat{{\rm Mat}}
\def \N{\mathbb N}
\def \o{{\mathfrak o}}
\def \ord{{\rm ord}}
\def \O{{\cal O}}
\def \p{{{\mathfrak p}}}
\def \P{{\mathbb P}}
\def \ph{\varphi}
\def \prf{\noindent{\bf Proof: }}
\def \Per{{\rm Per}}
\def \q{{\mathfrak q}}
\def \qed{\ifmmode\eqno $\square$\else\noproof\vskip 12pt plus 3pt minus 9pt \fi}
 \def\noproof{{\unskip\nobreak\hfill\penalty50\hskip2em\hbox{}%
     \nobreak\hfill $\square$\parfillskip=0pt%
     \finalhyphendemerits=0\par}}
\def \Q{\mathbb Q}
\def \res{{\rm res}}
\def \R{{\mathbb R}}
\def \Re{{\rm Re \hspace{1pt}}}
\def \r{{\mathfrak r}}
\def \ra{\rightarrow}
\def \rank{{\rm rank}}
\def \supp{{\rm supp}}
\def \SL{{\rm SL}}
\def \Spin{{\rm Spin}}
\def \t{{{\mathfrak t}}}
\def \T{{\mathbb T}}
\def \tr{{\hspace{1pt}\rm tr\hspace{2pt}}}
\def \vol{{\rm vol}}
\def \z{\zeta}
\def \Z{\mathbb Z}
\def \={\ =\ }

\newtheorem{theorem}{Theorem}[section]
\newtheorem{conjecture}[theorem]{Conjecture}
\newtheorem{lemma}[theorem]{Lemma}
\newtheorem{corollary}[theorem]{Corollary}
\newtheorem{proposition}[theorem]{Proposition}

\newcommand{\norm}[1]{\left|\hspace{-1pt}\left| #1\right|\hspace{-1pt}\right|}
\renewcommand{\matrix}[4]{\left( \begin{array}{rr}#1 & #2 \\ #3 & #4 \end{array}
            \right)}
\renewcommand{\sp}[1]{\left\langle #1 \right\rangle}

\section*{Introduction}
Let $G$ be a semisimple Lie group and $\Ga\subset G$ an arithmetic
subgroup. For a finite dimensional representation $(\rho,E)$ of $G$ the
cohomology groups
$H^\bullet(\Ga,E)$ are related  to automorphic forms and have
for this reason been studied by many authors.
The case of infinite dimensional representations has only very recently come into focus, mostly in connection with the Patterson Conjecture on the divisor of the Selberg zeta function \cite{BO1,BO2,BO3,D1,Juhl}. 
In this paper we want to show that the Patterson conjecture
\cite{BO1} is related to the Lewis correspondence \cite{LZ}, i.e.,
that the multiplicities of automorphic representations can be expressed in 
terms of cohomology groups with certain infinite dimensional coefficient
spaces.

One way to put (a special case of ) the Patterson conjecture 
for \emph{cocompact
torsion-free}
$\Ga$ in a split group $G$ is to say that the multiplicity $N_\Ga(\pi)$ of an
irreducible unitary principal series representation $\pi$ in the space
$L^2(\Ga\bs G)$ is given by
$$
N_\Ga(\pi)\=
\dim H^{d-r}(\Ga,\pi^\omega),
$$
where $r$ is the rank of $G$ and $\pi^\omega$ is the subspace of
analytic vectors in $\pi$, finally, $d=\dim(G/K)$ is the dimension of the symmetric space attached to $G$, where $K$ is a maximal compact subgroup.

Our main result states that this assertion can be generalized to all arithmetic
groups provided the ordinary group cohomology is replaced by the cuspidal
cohomology.
It will probably also work for more general lattices, but we stick to arihmetic groups, because some of the constructions used in this paper, like the Borel-Serre compactification, or the decomposition of the regular $G$-representation on the space $L^2(\Ga\bs G)$, have in the literature only been formulated for arithmetic groups.
The relation to the Lewis correspondence is as follows.
In \cite{Zagier} Don Zagier states that the correspondence for $\Ga={\rm
PSL}_2(\Z)$ can be interpreted as the identity
$$
N_{\Ga}(\pi)\=\dim H_{par}^1(\Ga,\pi^{\omega /2}),
$$
where $\pi$ is as before, $\pi^{\omega/2}$ is a slightly bigger space than
$\pi^\omega$ and $H_{par}^1$ is the parabolic cohomology. Since
$\pi$ is a unitary principal series representation it follows that
$N_\Ga(\pi)$ coincides with the multiplicity of $\pi$  in the cuspidal part
$L\cusp^2(\Ga\bs G)$ of
$L^2(\Ga\bs G)$. More precisely, the correspondence gives an isomorphism
$$
\Hom_G\left(\pi,L\cusp^2(\Ga\bs G)\right)\ra
H_{par}^1(\Ga,\pi^{\omega/2}).
$$
 As a consequence of our main result we will get the following theorem.

\begin{theorem}\label{main-Fuchs}
For every Fuchsian group $\Ga$ we have
$$
N_\Ga(\pi)\=\dim H\cusp^1(\Ga,\pi_{it}^{\omega}).
$$
\end{theorem}

Here $H\cusp^\bullet$ is the cuspidal cohomology. For finite dimensional
modules the cuspidal cohomology is a subspace of the parabolic cohomology. 

The following is our main theorem.

\begin{theorem}\label{main}
Let $\Ga$ be a
torsion-free arithmetic subgroup of a split semisimple Lie group $G$. Let $\pi\in\hat G$ be an
irreducible unitary principal series representation. Then
$$
N_{\Ga}(\pi)\= \dim H\cusp^{d-r}(\Ga,\pi^\omega),
$$
where $d=\dim G/K$ and $r$ is the real
rank of $G$.

For $G$ non-split the assertion remains true for a generic set of representations $\pi$.
\end{theorem}

This raises many questions. For a finite dimensional representation $E$ it is
known that the cuspidal cohomology is a subspace of the parabolic
cohomology. 
The same assertion for infinite dimensional $E$ is wrong in general, see Corollary \ref{4.2}.
Can one characterize those infinite dimensional $E$ for which the cuspidal cohomology indeed injects into ordinary cohomology?

Another question suggests itself:  in which sense does our construction in
the case
${\rm PSL}_2(\Z)$ coincide with the Lewis correspondence?
To even formulate a conjecture we must assume two further conjectures.
First assume that in the relevant cases
 cuspidal and parabolic cohomology coincide; next assume that the
parabolic cohomology with coefficients in
$\pi^{\omega}$ agrees with parabolic cohomology in $\pi^{\omega/2}$. Let
$M_\la$ be the space of cusp forms of eigenvalue $\la$. Then our
construction gives a map into the dual space of the cohomology,
$\alpha: M_\la\to H_{par}^1(\Ga,\pi^\omega)^*$. The Lewis
construction on the other hand gives a map $\beta: M_\la\to
H_{par}^1(\Ga,\pi^\omega)$. Together they define a duality on $M_\la$. One
is tempted to speculate that this duality coincides with the natural
duality given by the integral on the upper half plane. If that were so,
then the two maps $\alpha$ and $\beta$ would determine each other.

\section{Fuchsian groups}
Let $G$ be the group $\SL_2(\R)/\pm 1$.
For
$s\in\C$ let
$\pi_s$ denote the principal series representation with parameter $s$.
Recall that this representation can be viewed as the regular
representation on the space of square integrable sections of a line bundle
over $\P^1(\R)\cong G/P$, where $P$ is the subgroup of upper triangular
matrices. For
$s\in i\R$ this representation will be irreducible unitary. For any
admissible representation $\pi$ of $G$ let
$\pi^\omega$ denote the space of analytic vectors in $\pi$. Then
$\pi^\omega$ is a locally convex vector space with continuous
$G$-representation (\cite{KS}, p. 463). Let $\pi^{-\omega}$ be its
continuous dual. For $\pi=\pi_s$ the space $\pi_{s}^\omega$ is the
space of analytic sections of a line bundle over $\P^1(\R)$. Let
$\pi_{s}^{\omega/2}$ denote the space of sections which are smooth everywhere
and analytic up to the possible exception of finitely many points. Let
$\Ga=\SL_2(\Z)/\pm 1$ be the modular group. For an irreducible representation
$\pi$ of $G$ let
$N_\Ga(\pi)$ be its multiplicity in
$L^2(\Ga\bs G)$. Let
$H_{par}^1(\Ga,\pi_s^\omega)$ denote the \emph{parabolic cohomology}, i.e., the
subspace of $H^1(\Ga,\pi_s^\omega)$ generated by all cocycles $\mu$ which
vanish on parabolic elements. For the group $\Ga=\SL_2(\Z)/\pm 1$ this
means that $H_{par}^1$ consists of all cohomology classes which have a
representing cocycle $\mu$ with

$\mu\matrix 1101 =0$.
In \cite{Zagier} D. Zagier stated that for $s\in
i\R$,
$$
N_\Ga(\pi)\=\dim H_{par}^1(\Ga,\pi_{s}^{\omega/2}).
$$

We will first relate this to the Patterson Conjecture for the cocompact
case.

\begin{theorem}\label{2.1}
Let $\Ga\subset G$ be a discrete, cocompact and torsion-free subgroup,
 then for $s\in i\R$,
\begin{eqnarray*}
N_{\Ga}(\pi_{s})&=&\dim H_{par}^1(\Ga,\pi_{s}^\omega)\=\dim
H^1(\Ga,\pi_{s}^\omega).
\end{eqnarray*}
\end{theorem}

\prf
Since $\Ga$ does not contain parabolic elements the parabolic cohomology
coincides with the ordinary group cohomology. The Patterson Conjecture
\cite{BO1, D2} shows that
$$
N_\Ga(\pi_s)\=\dim
H^1(\Ga,\pi_{s}^{-\omega})-2\dim H^2(\Ga,\pi_{s}^{-\omega}).
$$
Poincar\'e duality \cite{BO2} implies that the dimension of
the space $H^j(\Ga,\pi_{s}^\omega)$ equals the dimension
of $H^{2-j}(\Ga,\pi_{s}^{-\omega})$. The Theorem follows from the next
lemma.

\begin{lemma}\label{vanish}
For every Fuchsian group we have $H^0(\Ga,\pi_{s}^\omega)=0$.
\end{lemma}

\prf
 For this recall that every $f\in\pi_s^\omega$
is a continuous function on $G$ satisfying among other things, $f(nx)=f(x)$ for
every $n\in N$, where $N$ is the unipotent group of all matrices modulo $\pm
1$ which are upper triangular with ones on the diagonal.
 If $f$ is $\Ga$-invariant,
then $f\in C(G/\Ga)$. By Moore's Theorem (\cite{Zimmer}, Thm. 2.2.6) it follows
that the action of $N$ on $G/\Ga$ is ergodic. In particular, this implies
that $f$ must be constant. Since $s\in i\R$ this implies that $f=0$.
\qed

\section{Arbitrary arithmetic groups}
Throughout, let $G$ be a semisimple Lie group with finite center and
finitely many connected components.

Let $\Ga$ be an arithmetic subgroup of  $G$ and assume that $\Ga$ is torsion-free. Then
$\Ga$ is the fundamental group of $\Ga\bs X$, where
$X=G/K$ the symmetric space and every $\Ga$-module $M$
induces a local system or locally constant sheaf $\CM$ on
$\Ga\bs X$. In the \'etale picture the sheaf $\CM$ equals $\CM=\Ga\bs
(X\times M)$, (diagonal action). Let
$\overline{\Ga\bs X}$ denote the Borel-Serre compactification \cite{BS} of
$\Ga\bs X$, then $\Ga$ also is the fundamental group of $\overline{\Ga\bs
X}$ and
$M$ induces a sheaf also denoted by $\CM$ on $\overline{\Ga\bs
X}$. This notation is consistent as the sheaf on ${\Ga\bs
X}$ is indeed the restriction of the one on
$\overline{\Ga\bs X}$. Let $\partial(\Ga\bs X)$ denote the
boundary of the Borel-Serre compactification. We have
natural identifications
$$
H^j(\Ga,M)\ \cong\ H^j(\Ga\bs X,\CM)\ \cong
H^j(\overline{\Ga\bs X},\CM).
$$
 We define the \emph{parabolic cohomology} of a $\Ga$-module $M$ to be the
kernel of the restriction to the boundary, ie,
$$
H_{par}^j(\Ga,M)\df \ker\left( H^j(\overline{\Ga\bs
X},\CM)\ra H^j(\partial (\Ga\bs X),\CM)\right).
$$
The long exact sequence of the pair $(\overline{\Ga\bs X},\partial(\Ga\bs X))$
gives rise to
$$
\dots\ra H_c^j(\Ga\bs X,\CM)\ra H^j(\Ga\bs X,\CM)=
$$ $$
=H^j(\overline{\Ga\bs X},\CM)\ra H^j(\partial(\Ga\bs
X),\CM)\ra\dots
$$
The image of the cohomology with compact supports under
the natural map is called the \emph{interior cohomology}
of $\Ga\bs X$ and is denoted by $H_!^j(\Ga\bs X,\CM)$. The
exactness of the above sequence shows that
$$
H_{par}^j(\Ga,M)\ \cong\ H_!^j(\Ga\bs X,\CM).
$$
Let $E$ be a locally convex space. We shall write $E'$ for its
topological dual. We assume that $\Ga$ acts linearly and continuously on $E$.
We will present a natural complex that computes the cohomology
$H^\bullet(\Ga,E)$. 

Let $\CE_\Ga$ be the locally constant sheaf on $\Ga\bs X$ given by $E$.
Then $\CE_\Ga$ has stalk
$E$ and
$H^\bullet(\Ga,E)=H^\bullet(X_\Ga,\CE_\Ga )$.

Let $\Omega_\Ga^0,\dots,\Omega_\Ga^d$ be the sheaves of
differential forms on $X_\Ga$ and let $\CE_\Ga^p$ be the
sheaf locally given by
$$
\CE_\Ga^p(U)\= \Omega_\Ga^p(U)\hat\otimes \CE_\Ga(U),
$$
where $\hat\otimes$ denotes the completion of the algebraic tensor product
$\otimes$ in the projective topology. 
Write $X_\Ga=\Ga\bs G/K=\Ga\bs X$.
Let $d$ denote the exterior
differential. Then
$D=d\otimes 1$ is a differential on $\CE_\Ga^\bullet$ and
$$
0\ra\CE_\Ga\stackrel D\ra\CE_\Ga^0\stackrel
D\ra\cdots\stackrel D\ra\CE_\Ga^d\ra 0
$$
is a fine resolution of $\CE_\Ga$. Hence
$H^\bullet(X_\Ga,\CE_\Ga)=H^\bullet(X_\Ga,\CE_\Ga^\bullet)$.

Let $\Omega^\bullet(X)$ be the space of differential forms
on $X$. The complex $\CE_\Ga^\bullet(X_\Ga)$ is isomorphic
to the space of $\Ga$-invariants
$(\Omega^\bullet(X)\hat\otimes E)^\Ga$. So we get
$$
H^\bullet(\Ga,E)\cong
H^\bullet\left((\Omega^\bullet(X)\hat\otimes E)^\Ga\right).
$$
We can write
$$
\Omega^p(X)\= (C^\infty(G)\otimes\wedge^p\p^*)^K,
$$
where $\p$ is the positive part in the Cartan
decomposition $\g=\k\oplus\p$, where $\g$ is the complexified Lie algebra of
$G$ and $\k$ is the complexified Lie algebra of $K$. The group
$K$ acts on
$\p^*$ via the coadjoint representation and on $C^\infty(G)$ via right
translations and
$\Ga$, or more precisely $G$, acts by left translations on $C^\infty(G)$.

From now on we assume that $E$ is not only a $\Ga$-module
but is a topological vector space that carries a continuous
$G$-representation. We say that
$E$ is
\emph{admissible} if every $K$-isotype $E(\tau)$, $\tau\in\hat K$ is finite
dimensional. Let
$E^\infty$ denote the subspace of smooth vectors. We say
that $E$ is smooth if $E=E^\infty$. We then have
$$
(C^\infty(G)\otimes\wedge^p\p^*)\hat\otimes E \cong
C^\infty(G)\hat\otimes (E\otimes\wedge^p\p^*)
$$
as a $G\times
K$-module, where $G$ acts diagonally  on $C^\infty(G)$ by left translations
and on
$E$ by the given representation. The group $K$ acts diagonally on
$C^\infty(G)$ by right translations and on $\wedge^p\p^*$ via the coadjoint
action.

\begin{lemma}
For any locally convex complete topological vector space $F$ we have
$$
C^\infty(G)\hat\otimes F\ \cong\ C^\infty(G,F),
$$
where the right hand side denotes the space of all smooth
maps from $G$ to $F$.
\end{lemma}

\prf See \cite{Gro}, Example 1 after Theorem 13.
\qed

Thus we have a $G\times K$-action on the space
$C^\infty(G,\wedge^p\p^*\otimes E)$ given by
$$
(g,k).f\= (\Ad^*(k)\otimes g)\, L_g\, R_k\, f,
$$
where $L_g f(x)=f(g^{-1}x)$ and $R_kf(x)=f(xk)$.

The map
$$
\psi\colon C^\infty(G,\wedge^p\p^*\otimes E) \ \ra\
C^\infty(G,\wedge^p\p^*\otimes E)
$$
given by
$$
\psi(f)(x)=(1\otimes x^{-1}).f(x)
$$

is an isomorphism to the same space with a different the
$G\times K$ structure. 
Indeed, one computes,
\begin{eqnarray*}
\psi((g,k).f)(x) &=& (1\otimes x^{-1})(g,k).f(x)\\
&=& (1\otimes x^{-1}) (\Ad^*(k)\otimes g)f(g^{-1}xk)\\
&=& (\Ad^*(k)\otimes k(g^{-1}xk)^{-1}) f(g^{-1}xk)\\
&=& (\Ad^*(k)\otimes k) R_k L_g \psi(f)(x).
\end{eqnarray*}

For a smooth $G$-representation $F$ we write
$H^\bullet(\g,K,F)$ for the cohomology of the standard complex of
$(\g,K)$-cohomology \cite{BorWall}. Then $H^\bullet(\g,K,F)=
H^\bullet(\g,K,F_K)$, where $F_K$ is the $(\g,K)$-module of $K$-finite
vectors in $F$.

If we assume that $E$ is smooth, we
get from this
$$
H^\bullet(\Ga,E) \= H^\bullet(\g,K,C^\infty(\Ga\bs
G)\hat\otimes E).
$$
In the case of finite dimensional $E$ one can replace $C^\infty(\Ga\bs G)$ with the space of functions of moderate growth \cite{Borel-stableII}.
This is of importance, since it leads to a decomposition of the cohomology space into the cuspidal part and contributions from the parabolic subgroups.
To prove this, one starts with differential forms of moderate growth and applies $\psi$.
For infinite dimensional  $E$ this proof does not work, since it is not clear that $\psi$ should preserve moderate growth, even if one knows that the matrix coefficients of $E$ have moderate growth.

By the Sobolev Lemma the space of smooth vectors $L^2(\Ga\bs G)^\infty$ of the
natural unitary representation of $G$ on $L^2(\Ga\bs G)$ is a subspace of
$C^\infty(\Ga\bs G)$. The representation $L^2(\Ga\bs G)$ splits as $L^2(\Ga\bs
G)=L_{disc}^2\oplus L_{cont}^2$, where $L_{disc}^2=\bigoplus_{\pi\in\hat G}
N_\Ga(\pi)\pi$ is a direct Hilbert sum of irreducible representations and
$L_{cont}^2$ is a finite sum of continuous Hilbert integrals. The space of cusp
forms $L\cusp^2(\Ga\bs G)=\bigoplus_{\pi\in\hat G} N_{\Ga,\rm
cusp}(\pi)\pi$ is a subspace of $L_{disc}^2$. Note that $L\cusp^2(\Ga\bs
G)^\infty$ is a closed subspace of $C^\infty(G)$. The \emph{cuspidal
cohomology} is defined by
$$
H\cusp^\bullet(\Ga,E)\ =\ H^\bullet(\g,K,L\cusp^2(\Ga\bs
G)^\infty\hat\otimes E).
$$
For finite dimensional $E$ it turns out that $H\cusp^\bullet(\Ga,E)$
coincides with the image in $H^\bullet(\g,K,C^\infty(\Ga\bs G)\hat\otimes
E)$ under the inclusion map. This comes about as a consequence of the fact
that the cohomology can also be computed using functions of uniform
moderate growth and that in the space of such functions, $L^2\cusp(\Ga\bs
G)^\infty$ has a $G$-complement. The Borel-conjecture \cite{F} is a
refinement of this assertion.
For infinite dimensional $E$ this injectivity does not hold in general, see Corollary \ref{4.2}.

We define the \emph{reduced cuspidal cohomology} to be the image $\tilde
H^\bullet\cusp(\Ga,E)$ of $H^\bullet\cusp(\Ga,E)$ in
$H^\bullet(\Ga,E)$. Finally, let $H_{(2)}(\Ga,E)$ be the image of the space
 $H^\bullet(\g,K,L^2(\Ga\bs G)^\infty\hat\otimes E)$ in
$H^\bullet(\g,K,C^\infty(\Ga\bs G)\hat\otimes E)$. 

\begin{proposition}
We have the following inclusions of cohomology groups,
$$
\tilde H\cusp^\bullet(\Ga,E)\ \subset\ H_{par}^\bullet(\Ga,E)\ \subset\
H_{(2)}^\bullet(\Ga,E).
$$
\end{proposition}

\prf
The cuspidal condition ensures that every cuspidal class vanishes on each
homology class of the boundary. This implies the first conclusion. Since
every parabolic class has a compactly supported representative, the second
also follows.
\qed

\section{Gelfand Duality}
Recall that a \emph{Harish-Chandra module} is a $(\g,K)$-module which is
admissible and finitely generated. Every Harish-Chandra module is of
finite length. For a Harish-Chandra module
$V$ write
$\tilde  V$ for its dual, ie, $\tilde V=(V^*)_K$, the $K$-finite vectors
in the algebraic dual. 

A \emph{globalization} of a Harish-Chandra module $V$ is a continuous
representation of $G$ on a complete locally convex vector space $W$ such that
$V$ is isomorphic to the $(\g,K)$-module of $K$-finite vectors $W_K$. It
was shown in \cite{KS} that there is a minimal globalization $V^{\rm min}$ and
a maximal globalization $V^{\rm max}$ such that for every globalization $W$
there are unique functorial continuous linear $G$-maps
$$
V^{\rm min}\ \hookrightarrow\ W\ \hookrightarrow\ V^{\rm max}.
$$
The spaces $V^{\rm min}$ and $V^{\rm max}$ are given explicitly by
$$
V^{\rm min}\= C_c^\infty(G) \otimes_{\g,K} V
$$
and
$$
V^{\rm max}\= \Hom_{\g,K}(\tilde V,C^\infty(G)).
$$
The action of $G$ on $V^{max}$ is given by
$$
g.\al(\tilde v)(x)\=\al(\tilde v)(g^{-1}x).
$$
Let $\hat G$ be the unitary dual of $G$, i.e., the set of all isomorphism
classes of irreducible unitary representations of $G$. Note \cite{KS} that for
$\pi\in\hat G$ we have
$(\pi_K)^{\rm min}=\pi^\omega$ and
$(\pi_K)^{\rm max}=\pi^{-\omega}$.

The following is a key result of this paper.

\begin{theorem}\label{3.7}
Let $F$ be  a smooth $G$-representation on a complete locally convex
topological vector space. Then there is a functorial isomorphism
$$
H^\bullet(\g,K,F\hat\otimes V^{\rm max})\ \ra\
\Ext_{\g,K}^\bullet(\tilde V,F),
$$
where as usual one writes $\Ext_{\g,K}^\bullet(\tilde V,F)$ for $\Ext_{\g,K}^\bullet(\tilde V,F_K)$.
\end{theorem}

\prf  We have
\begin{eqnarray*}
F\hat\otimes V^{\rm max} &=& F\hat\otimes \Hom_{\g,K}(\tilde V,C^\infty(G))\\
&=& \Hom_{\g,K}(\tilde V,F\hat\otimes C^\infty(G))\\
&=& \Hom_{\g,K}(\tilde V, C^\infty(G,F)).
\end{eqnarray*}

\begin{lemma}
The map
\begin{eqnarray*}
\Hom_{\g,K}\left(\tilde V,C^\infty(G,F))\right)^{\g,K}
&\to&\Hom_{\g,K}\left(\tilde V, F\right)\\
\phi &\mapsto &\al,
\end{eqnarray*}
with $\al(\tilde v)=\phi(\tilde v)(1)$ is an isomorphism.
\end{lemma}

\prf
Note that $\phi$ satisfies
$$
\phi(X.\tilde v)(x)\=\left.\frac d{dt}\phi(\tilde
v)(x\,\exp(tX))\right|_{t=0},\qquad X\in\g,
$$
since it is a $(\g,K)$-homomorphism. Further,
$$
\left.\frac d{dt} \phi(\tilde v)(\exp(tX)\, x)\right|_{t=0}\= X.\phi(\tilde
v)(x),\qquad X\in\g,
$$
since $\phi$ is $(\g,K)$-invariant. Similar identities hold for the $K$-action.
This implies that $\al$ is a $(\g,K)$-homomorphism. Note that the
$(\g,K)$-invariance of $\phi$ also leads to
$$
\phi(\tilde v)(gx)\= g.\phi(\tilde v)(x),\qquad g,x\in G.
$$
Hence if $\al =0$ then $\phi=0$ so the map is injective. For surjectivity let
$\al$ be given and define $\phi$ by $\phi(\tilde v)(x)=x.\al(\tilde v)$. Then
$\phi$ maps to $\al$.
\qed

By the Lemma we get an isomorphism
$$
H^0(\g,K,F\hat\otimes V^{\rm max})\ \cong\ \Hom_{\g,K}(\tilde
V, F)\ \cong\ \Ext_{\g,K}^0(\tilde
V, F)
$$
and thus a functorial isomorphism on the zeroth level. We will show that both
sides in the theorem define universal $\delta$-functors \cite{L}. From this the
theorem will follow.
Fix $V$ and let $S^j(F)=H^j(\g,K,F\hat\otimes V^{\rm max})$ as well as
$T^j(F)=\Ext_{\g,K}^j(\tilde V,F)$. We will show that $S^\bullet$ and
$T^\bullet$ are universal $\delta$-functors from the category $Rep_s^\infty(G)$
defined below to the category of complex vector spaces. The objects of
$Rep_s^\infty(G)$ are smooth continuous representations of $G$ on Hausdorff
locally convex topological vector spaces and the morphisms are \emph{strong
morphisms}. A continuous $G$-morphism $f:A\ra B$ is called strong morphism or
\emph{s-morphism} if (a) $\ker f$ and $\im f$ are closed topological direct
summands and (b) $f$ induces an isomorphism of $A/\ker f$ onto $f(A)$. Then by
\cite{BorWall}, Chapter IX, the category $Rep_s^\infty(G)$ is an abelian
category with enough injectives. In fact, for $F\in Rep_s^\infty(G)$ the map
$F\ra C^\infty(G,F)$ mapping $f$ to the function $\al(x)=x.f$ is a monomorphism
into the s-injective object $C^\infty(G,F)$ (cf. \cite{BorWall}, Lemma
IX.5.2), which  is considered a $G$-module via
$x\al(y)=\al(yx)$.

Let us consider $S^\bullet$ first. By Corollary IX.5.6 of \cite{BorWall} we have
$$
S^\bullet(F)\= H^\bullet(\g,K,F\hat\otimes V^{\rm max})\ \cong\
H_d^\bullet(G,F\hat\otimes V^{\rm max}),
$$
where the right hand side is the differentiable cohomology. The functor
$.\hat\otimes V^{\rm max}$ is s-exact and therefore $S^\bullet$ is a
$\delta$-functor.
We show that it is erasable. For this it suffices to show that
$S^j(C^\infty(G,F))=0$ for $j>0$.
Now
$$
C^\infty(G,F)\hat\otimes V^{\rm max}\ \cong\ C^\infty(G)\hat\otimes F\hat\otimes
V^{\rm max}\ \cong\ C^\infty(G,F\hat\otimes V^{\rm max})
$$
and therefore for $j>0$,
$$
S^j(C^\infty(G,F))\ \cong\ H_d^j(G,C^\infty(G,F\hat\otimes V^{\rm max}))\= 0,
$$
since $C^\infty(G,F\hat\otimes V^{\rm max})$ is s-injective. Thus $S^\bullet$ is
erasable and therefore universal.

Next consider $T^\bullet(F)=\Ext_{\g,K}^\bullet(\tilde V,F)$. Since an exact
sequence of smooth representations gives an exact sequence of
$(\g,K)$-modules, it follows that $T^\bullet$ is a $\delta$-functor.To show
that it is erasable let $j>0$. Then
\begin{eqnarray*}
T^j(C^\infty(G,F)) &=& \Ext_{\g,K}^j(\tilde V, C^\infty(G)\hat\otimes F)\\
&=& H^j(\g,K,\Hom_\C(\tilde V,C^\infty(G))\hat\otimes F)\\
&=& H_d^j(G,\Hom_\C(\tilde V,C^\infty(G))\hat\otimes F)\\
&=& H_d^j(G,\Hom_\C(\tilde V,C^\infty(G)))\hat\otimes F\\
&=& \Ext_{\g,K}^j(\tilde V,C^\infty(G))\hat\otimes F.
\end{eqnarray*}
By Theorem 6.13 of \cite{KS} we have $\Ext_{\g,K}^j(\tilde V,C^\infty(G))=0
$.
The Theorem is proven.
\qed

Choosing  $C^\infty(\Ga\bs G)$ and $L^2_{\mathrm{cusp}}(\Ga\bs G)^\infty$ for
$F$ in  Theorem \ref{3.7} gives the following
Corollary. 
\begin{corollary}\label{GelfandDuality}
\begin{enumerate}
\item[{\rm(i)}]
$$
H^p(\Ga, V^{\rm max})\ \cong\ \Ext_{\g,K}^p(\tilde V,C^\infty(\Ga\bs
G)).
$$
For $\Ga$ cocompact and $p=0$ this is known under the name \emph{Gelfand
Duality}.
\item[{\rm(ii)}]
$$
H\cusp^\bullet(\Ga,V^{\rm max})\ \cong\
\Ext_{\g,K}^\bullet(\tilde V,L\cusp^2(\Ga\bs G)^\infty).
$$
\end{enumerate}
\end{corollary}

\section{The case of the maximal globalization}
The space of cusp forms decomposes discretely,
$$
L\cusp^2(\Ga\bs G)\=\bigoplus_{\pi\in\hat G} N_{\Ga, \rm cusp}(\pi)\pi.
$$
Suppose that $V$ has an infinitesimal character $\chi$. Let $\hat G(\chi)$ be
the set of all irreducible unitary representations of $G$ with infinitesimal
character $\chi$. It is easy to see that
\begin{eqnarray*}
\Ext_{\g,K}^\bullet(\tilde V,L\cusp^2(\Ga\bs G)^\infty)
&=&\Ext_{\g,K}^\bullet\left(\tilde
V,\bigoplus_{\pi\in \hat G(\chi)} N_{\Ga,\rm cusp}(\pi)\pi_K\right)\\
&=&\bigoplus_{\pi\in \hat G(\chi)} N_{\Ga,\rm cusp}(\pi)\, \Ext_{\g,K}^\bullet\left(\tilde
V,\pi_K\right)\\
&=&\bigoplus_{\pi\in \hat G(\chi)} N_{\Ga,\rm cusp}(\pi)\, \Ext_{\g,K}^\bullet\left(\tilde
\pi_K,V\right)
\end{eqnarray*}
The last line follows by dualizing.

Let $P$ be a parabolic subgroup and $\m\oplus\a\oplus \n $ a Langlands
decomposition of its Lie algebra.
\begin{lemma}\label{Frobenius}
For a $(\g,K)$-module $\pi$ and a $(\a\oplus \m,K_M)$-module $U$ we have
$$
\Hom_{\g,K}(\pi,\Ind_P^G(U))\ \cong\
\Hom_{\a\oplus\m,K_M}(H_0(\n,\pi),U\otimes
\C_{\rho_P}),
$$
where $\C_{\rho_P}$ is the one dimensional $A$-module given by $\rho_P$.
\end{lemma}

\prf See \cite{HeSch} page 101.
\qed

\begin{lemma}\label{Spektralsequenz}
Let $\CC$ be an abelian category with enough injectives. Let $\a$ be a finite
dimensional abelian complex Lie algebra and let $T$ be a covariant left exact
functor from $\CC$ to the category of $\a$-modules. Assume that $T$ maps
injectives to $\a$-acyclics and that $T$ has finite cohomological
dimension, i.e., that $R^pT=0$ for $p$ large.
Let $M$ be an object of $\CC$ such that $R^pT(M)$ is finite dimensional for every
$p$. Let $H_\a$ denote the functor $H^0(\a,\cdot)$. Then
$$
\sum_{p\ge 0}\binom pr (-1)^{p+r}\,\dim R^p(H_\a\circ T)(M)\=\sum_{p\ge
0}(-1)^p\,\dim H^0(\a,R^pT(M)),
$$
where $r=\dim\a$.
If $M$ is $T$-acyclic, then these alternating sums degenerate to
$$
\dim R^r(H_\a \circ T)(M)\= \dim H^0(\a,T(M)).
$$
\end{lemma}

\prf
Split $\a=\a_1\oplus \b_1$, where $\dim\a_1=1$. Consider the Grothendieck spectral
sequence with $E_2^{p,q}=H^p(\a_1,R^q(H_{\b_1}\circ T)(M))$ which abuts to
$R^{p+q}(H_\a\circ T)(M)$ (see \cite{L}, Theorem XX.9.6).
Since $\a_1$ is one dimensional, for any finite dimensional $\a_1$-module $V$ we
have $H^0(\a,V)\cong H^1(\a,V)$ and $H^p(\a_1,V)=0$ if $p>1$. This implies
$E_2^{0,q}\cong E_2^{1,q}$ and $E_2^{p,q}=0$ for $p\notin\{ 0,1\}$. The spectral
sequence therefore degenerates and
\begin{eqnarray*}
&&\sum_{p\ge 0}\binom pr (-1)^{p+r}\dim R^p(H_\a\circ T)(M)\\
&=& \sum_{p\ge 0}\binom pr (-1)^{p+r}\dim E_2^{0,p}+
    \sum_{p\ge 1}\binom pr (-1)^{p+r} \dim E_2^{1,p-1}\\
&=& \sum_{p\ge 0}\binom pr (-1)^{p+r}\dim E_2^{0,p}+
    \sum_{p\ge 0}\binom {p+1}r (-1)^{p+r-1} \dim E_2^{1,p}\\
&=& \sum_{p\ge 0}\left(\binom {p+1}r -\binom pr\right) (-1)^{p+r-1}\dim E_2^{0,p}\\
&=& \sum_{p\ge 0}\binom {p}{r-1} (-1)^{p+r-1}\dim E_2^{0,p}\\
&=& \sum_{p\ge 0}\binom p{r-1}(-1)^{p+r-1}\dim H^0(\a_1,R^p(H_{\b_1}\circ T)(M)).
\end{eqnarray*}
Next we split $\b_1=\a_2\oplus\b_2$, where $\a_2$ is one-dimensional. Since the
$\a_1$-action commutes with the $a_2$-action the isomorphism
$$
H^0(\a_2,R^p(H_{\b_2}\circ T)(M))\ \cong\ H^1(\a_2,R^p(H_{\b_2}\circ T)(M))
$$
is an $\a_1$-isomorphism. Therefore we apply the same argument to get down to
$$
\sum_{p\ge 0} \binom p{r-2} (-1)^{p+r-2}\dim H^0(\a_1\oplus\a_2,R^p(H_{\b_2}\circ
T)(M)).
$$
Iteration gives the claim.

To get the last assertion of the lemma, note that if $M$ is $T$-acyclic, then following the inductive argument above, one sees that $R^p(H_\a \circ T(M)=0$ for $p>r$.
\qed

Let $P$ be a minimal parabolic subgroup of $G$ so that $M$ is compact.
For a unitary irreducible representation $\sigma$ of $M$ and linear functional 
$\nu\in i\a^*$ we obtain the unitary principal series representation
$\pi_{\sigma,\nu}$ of $G$ induced from $P$.
Let $\t$ be a Cartan subalgebra of $\m={\rm Lie}_\C(M)$.
Then $\h=\a\oplus\t$ is a Cartan subalgebra of $\g$.
Let $\La_\sigma\in\t^*$ be a representative of the infinitesimal character of $\sigma$.
Then $\La_\sigma +\nu\in\h^*$ is a representative of the infinitesimal character of $\pi_{\sigma,\nu}$.
We say that the parameters $(\sigma,\nu)$ are \emph{generic} if $\pi_{\sigma,\nu}$ is irreducible and for any two $w,w'$ in the Weyl group of $(\h,\g)$ the linear functional 
$$
w(\La_\sigma +\nu)|_\a -w'(\La_\sigma +\nu)|_\a
$$
on $\a$ is not a positive integer linear combination of positive roots. 

\begin{theorem}
If $G$ is $\R$-split and $\pi_{\sigma,\nu}$ is irreducible, we have
$$
N_{\Ga}(\pi_{\sigma,\nu})\= H\cusp^r(\Ga,\pi_{\sigma,\nu}^{-\omega}).
$$
If $G$ is not split, the assertions remains true if the parameters $(\sigma,\nu)$ are generic.
\end{theorem} 

\prf
First note that since $G$ is split, the decomposition of $L^2(\Ga\bs G)$ as in
\cite{Langlands,Arthur} implies that for imaginary $\nu$ one has
$N_\Ga(\pi_{\sigma,\nu})=N_{\Ga,\rm cusp}(\pi_{\sigma,\nu})$, since the
Eisenstein series are regular at imaginary $\nu$. Applying Lemma
\ref{Frobenius} with
$\pi=L\cusp^2(\Ga\bs G)^\infty$ and $\Ind_P^G(U)=\pi_{\sigma,\nu}$ we find
\begin{eqnarray*}
&&\hspace{-40pt}\Hom_{\a\oplus\m,M}\left(H_0(\n,L\cusp^2(\Ga\bs
G)^\infty_K),\sigma\otimes(\nu+\rho_P)\right)\\ 
&\cong&\Hom_{\g,K}\left(L\cusp^2(\Ga\bs
G)^\infty,\pi_{\sigma,\nu}\right)\\ &\cong&\Hom_{\g,K}\left(\tilde
\pi_{\sigma,\nu},L\cusp^2(\Ga\bs G)^\infty\right)
\end{eqnarray*}
so that
$$N_{\Ga, \cusp}(\pi_{\sigma,\nu})=
\dim \Hom_{\a\oplus\m,M}\left(H_0(\n,L\cusp^2(\Ga\bs G)^\infty_K),
\sigma\otimes(\nu+\rho_P)\right).$$
In order to calculate the latter
we apply Lemma \ref{Spektralsequenz} to the category $\CC$ of
$(\g,K)$ modules and 
\begin{eqnarray*}
T(V) &=& \Hom_M(H_0(\n,\tilde V),U\otimes\C_{\rho_P})\\
&=&
(H_0(\n,\tilde V)^*\otimes U\otimes\C_{\rho_P})^M.
\end{eqnarray*}
The conditions of the Lemma \ref{Spektralsequenz}
are easily seen to be satisfied since $H_0(\n,\cdot)$ maps injectives to
injectives and $H^0(M,\cdot)$ is exact. Note that in the case of a
representation
$\pi$ of
$G$,
\begin{eqnarray*}
H_\a\circ T(\pi) &\cong&
\Hom_{\a\oplus\m,M}(H_0(\n,\tilde\pi_K),U\otimes\C_{\rho_P})\\ 
&\cong&
\Hom_{\g,K}(\tilde\pi,\Ind_P^G(U))
\end{eqnarray*}
by Lemma \ref{Frobenius}. From this we obtain
$$\Ext_{\g,K}^j(\tilde\pi,\Ind_P^G(U))=R^p(H_\a\circ T)(\pi).$$

Now Lemma \ref{Spektralsequenz} shows that
$$
\dim \Ext_{\g,K}^r(\tilde\pi,\Ind_P^G(U))
$$
equals
$$
\dim \Hom_{AM}(H_0(\n,\tilde\pi_K),U\otimes \C_{\rho_P}).
$$
Suppose that $\pi$ is an irreducible summand in $L\cusp^2(\Ga\bs G)$.
If $G$ is split, the set of $\pi$ which share the same infinitesimal
character as $\pi_{\sigma,\nu}$ equals the set of all $\pi_{\xi,w\nu}$, 
where $w$ ranges over the Weyl group and $\xi\in\hat M$.
Then the space
$\Hom_{AM}(H_j(\n,\tilde\pi_K),U\otimes \C_{\rho_P})$ is only non-zero for $\pi=\pi_{\xi,w\nu}$ for some $w$. But then Proposition 2.32 of \cite{HeSch} implies
that $H_j(\n,\tilde\pi_K)$ is zero unless $j=0$.
The same conclusion is assured in the non-split case by the genericity condition.
Now the proof is completed by the following calculation:
\begin{eqnarray*}
N_{\Ga, \mathrm{cusp}}(\pi_{\sigma,\nu})
&=&\dim \Hom_{\a\oplus\m,M}(H_0(\n,L\cusp^2(\Ga\bs G)^\infty_K),
   \sigma\otimes(\nu+\rho_P))\\
&=&\dim\Ext_{\g,K}^r(L\cusp^2(\Ga\bs G),\pi_{\sigma,\nu})\\
&=& \dim \Ext_{\g,K}^r(\tilde\pi_{\sigma,\nu},L\cusp^2(\Ga\bs G))\\
&=& \dim H\cusp^r(\Ga,\pi_{\sigma,\nu}^{-\omega}),
\end{eqnarray*}
where the last equality is a consequence of Corollary \ref{GelfandDuality}(ii),
applied to $V=\pi_{\sigma,\nu}$.
\qed

\section{Poincar\'e duality}
In order to conclude the main Theorem it suffices to prove the following
Poincar\'e duality.

\begin{theorem}\label{3.12}
(Poincar\'e duality)\\
For every Harish-Chandra module $V$,
$$
H\cusp^j(\Ga,V^{\rm max})\ \cong\
H\cusp^{d-j}(\Ga,\tilde V^{\rm min})^*,
$$
and both spaces are finite dimensional.
\end{theorem}

Before we prove the theorem, we add a Corollary.

\begin{corollary}
\label{4.2}
Let $\Ga$ be a torsion-free non-uniform lattice in $G={\rm PSL}_2(\R)$ and $\pi\in\hat G$ a principal series representation  with $N_{\Ga,{\rm cusp}}(\pi)\ne 0$.
Let $E=\tilde \pi_K^{min}$. Then the natural map $H\cusp^\bullet(\Ga,E)\to H^\bullet(\Ga,E)$ is not injective.
\end{corollary}

{\bf Proof of the Corollary:}
We have $H_{\rm cusp}^0(\Ga, \pi^{max})\ne 0$ and by the Poincar\'e duality, $H_{\rm cusp}^2(\Ga,E)\ne 0$.
However, as the cohomological dimension of $\Ga$ is $1$, it follows that $H^2(\Ga,E)=0$.
\qed

{\bf Proof of the Theorem:}
A \emph{duality} between two complex vector spaces $E,F$ is a bilinear
pairing,
$$
\sp{.,.}\colon E\times F\ \ra\ \C
$$
which is \emph{non-degenerate}, i.e., for every $e\in E$ and every
$f\in F$,
\begin{eqnarray*}
\sp{e,F}=0 &\Rightarrow & e=0,\\
\sp{E,f}=0 &\Rightarrow & f=0.
\end{eqnarray*}
We say that $E$ and $F$ are \emph{in duality} if there is a duality between
them.
Note that if $E$ and $F$ are in duality and one of them is finite dimensional,
then the other also is and their dimensions agree.
The pairing is called \emph{perfect} if it induces isomorphisms $E\cong F^*$
and $F\cong E^*$. If $E$ and $F$ are topological vector spaces then the pairing
is called \emph{topologically perfect} if it induces
topological isomorphisms
$E\cong F'$ and $F\cong E'$, where the dual spaces are equipped with the
strong dual topology. 

Now suppose that $V$ and $W$ are $\g,K$-modules in duality through a
$\g,K$-invariant pairing. Recall the canonical complex defining
$\g,K$-cohomology which is given by
$C^q(V)=\Hom_K(\wedge^q(\g/\k),V)=\left(\wedge^q(\g/\k)^*\otimes V\right)^K$.
Let $d=\dim G/K$. The prescription
$\sp{y\otimes v, y'\otimes w}=(-1)^q\sp{v,w}y\wedge y'$ gives a pairing
from
$C^q(V)\times C^{d-q}(W)$ to $\wedge^d(\g/\k)^*\cong\C$. Let $d:C^q\ra C^{q+1}$
be the differential, then one sees \cite{BorWall},
$\sp{da,b}=\sp{a,db}$.

Let $\pi$ be an irreducible unitary representation of $G$. Then the spaces
$\pi^{\infty}$ and
$\tilde\pi^{-\infty}$ are each other's strong duals \cite{C}. The same holds for
$V^{\rm max}$ and
$\tilde V^{\rm min}$
\cite{KS}. 

\begin{lemma}
The spaces $A=\pi^{-\infty}\hat\otimes V^{\rm max}$ and
$B=\tilde\pi^\infty\hat\otimes \tilde V^{\rm min}$ are each other's strong
duals. Both of them are LF-spaces.
\end{lemma}

\prf
Since $C^\infty(G)$ is nuclear and Fr\'echet and $\tilde \pi$ 
is a Hilbert space
the results of \cite{Tr}, \S III.50, allow us to conclude that
$C^\infty(G,\tilde \pi)'=C^\infty(G)\hat\otimes \tilde \pi$ is nuclear which 
is then 
true also for $C^\infty(G,\tilde \pi)$. Now the embedding of $\tilde \pi^\infty$
into $C^\infty(G,\tilde \pi)$ shows the nuclearity of $\tilde\pi^\infty$.

Since $\tilde V$ is finitely generated one can
embed the space 
$$
V^{\rm max}=\Hom_{\g,K}(\tilde V, C^\infty(G))
$$ 
into a strict
inductive  limit $\underset{\longrightarrow_j}{\lim}\ \Hom(\tilde
V^j,C^\infty(G))$ with finite  dimensional $V^j$'s. Then the nuclearity of
$V^{\rm max}$ follows from the nuclearity  of 
$$
\Hom(\tilde V^j,C^\infty(G))=(\tilde V^j)^*\otimes C^\infty(G).
$$

We conclude that the spaces $V^{\rm max}$ and $\tilde\pi^\infty$ are nuclear
Fr\'echet spaces.  Their
duals $\pi^{-\infty}$ and $\tilde V^{\rm min}$ are LF-spaces (see \cite{Gro},
Introduction IV). Therefore they all are barreled (\cite{Sch}, p. 61).  By
\cite{Sch}, p. 119 we know that the inductive completions of the tensor products
$\pi^{-\infty}\bar\otimes V^{\rm max}$ and
$\pi^\infty\bar\otimes V^{\rm min}$ are barreled. Since $V^{\rm max}$ and
$\pi^\infty$ are nuclear, these inductive completions coincide with the
projective completions. So $A$ and $B$ are barreled.
By Theorem 14 of \cite{Gro} it follows that the strong duals $A'$ and $B'$
are complete and by the Corollary to Lemma 9 of \cite{Gro} it follows that
$A'=B$ and $B'=A$. Finally, Lemma 9 of \cite{Gro} implies that $A$ and $B$
are LF-spaces.
\qed                                  

\begin{proposition}
For every $\pi\in\hat G$ and every Harish-Chandra module $V$ the vector spaces
$H^q(\g,K,\pi^{-\infty}\hat\otimes V^{\rm max})$ and
$H^q(\g,K,\tilde\pi^{\infty}\hat\otimes \tilde V^{\rm min})$ are finite
dimensional. The above pairing between their canonical complexes induces a
duality between them, so
$$
H^q(\g,K,\pi^{-\infty}\hat\otimes V^{\rm max})\ \cong\
H^{d-q}(\g,K,\tilde\pi^{\infty}\hat\otimes \tilde V^{\rm min})^*.
$$
\end{proposition}

\prf
Note that by Theorem \ref{3.7},
$$
H^q(\g,K,\pi^{-\infty}\hat\otimes V^{\rm max})\cong\Ext_{\g,K}^q(\tilde
V,\pi^{-\infty})\cong\Ext_{\g,K}^q(\tilde
V,\pi_K)
$$ 
and the latter space is finite dimensional (\cite{BorWall},
Proposition I.2.8). The proposition will thus follow from the next lemma.

\begin{lemma}
Let $A,B$ be smooth representations of $G$. Suppose that $A$ and $B$ are
LF-spaces and that they are in perfect topological
duality through a
$G$-invariant pairing. Assume that $H^\bullet(\g,K,A)$ is finite
dimensional. Then the natural pairing between $H^q(\g,K,A)$ and
$H^{d-q}(\g,K,B)$ is perfect.
\end{lemma}

\prf
We only have to show that
the pairing is non-degenerate. We will start by showing that the induced
map
$H^{d-q}(\g,K,B)$ to $H^q(\g,K,A)^*$ is injective. So let
$b\in Z^{d-q}(B)=C^{d-q}(B)\cap\ker d$ with
$\sp{a,b}=0$ for every
$a\in Z^{q}(A)$. Define a map $\psi\colon d(C^{q}(A))\ra \C$ by
$$
\psi(da)\= \sp{a,b}.
$$
We now show that the image $d(C^{q}(A))$ is a closed subspace of
$C^{q+1}(A)$ and that the map $C^{q}(A)/\ker d\ra d(C^{q}(A))$ is a
topological isomorphism. For this let $E$ be a finite dimensional subspace
of
$Z^{q+1}(A)$ that bijects to $H^{q+1}(\g,K,A)$. Since $E$ is finite
dimensional, it is closed. The map $\eta=d+1\colon C^{q}(A)\oplus E\ra
Z^{q+1}(A)$ is continuous and surjective. Since $C^{q}(A)$ and
$Z^{q+1}(A)$ are LF-spaces, the map $\eta$ is open (see \cite{Tr}, p. 78), hence
it induces a topological isomorphism
$(C^{q}(A)/\ker d)\oplus E\ra Z^{q+1}(A)$.
This implies that $d(C^{q}(A))$ is closed and $C^{q}(A)/\ker d\ra
d(C^{q}(A))$ is a topological isomorphism.
Consequently, the map $\psi$ is continuous. Hence it extends to a
continuous linear map on $C^{q+1}(A)$. Therefore, it is given by an
element $f$ of $C^{d-q-1}(A)$, so
$$
\sp{a,b}\= \sp{da,f}\= \sp{a,df}
$$
for every $a\in C^{q}(A)$. We conclude $b=df$ and thus the
non-degeneracy on one side. In particular it follows that
$H^\bullet(\g,K,B)$ is finite dimensional as well. The claim now follows
by symmetry.
\qed

We will now deduce Theorem \ref{3.12}. We have
\begin{eqnarray*}
H\cusp^q(\Ga,V^{\rm max}) &\cong & \bigoplus_{\pi\in\hat
G(\chi)}N_{\Ga,\rm cusp}(\pi)\,\Ext_{\g,K}^q(\tilde V,\pi_K)\\
&\cong & \bigoplus_{\pi\in\hat
G(\chi)}N_{\Ga,\rm cusp}(\pi)\,H^q(\g,K,\pi^{-\infty}\hat\otimes V^{\rm max})\\
&\cong & \bigoplus_{\pi\in\hat
G(\chi)}N_{\Ga,\rm cusp}(\pi)\,H^{d-q}(\g,K,\tilde\pi^{\infty}\hat\otimes
\tilde V^{\rm min})^*\\
&\cong & \bigoplus_{\pi\in\hat
G(\chi)}N_{\Ga,\rm cusp}(\pi)\,H^{d-q}(\g,K,\pi^{\infty}\hat\otimes
\tilde V^{\rm min})^*\\
&\cong& H^{d-q}\cusp(\Ga,\tilde V^{\rm min})^*.
\end{eqnarray*}
In the second to last step we have used the fact that $L\cusp^2$ is
self-dual. Theorem \ref{3.12} and thus Theorem \ref{main} follow.

It remains to deduce Theorem \ref{main-Fuchs}. For $\Ga$ torsion-free
arithmetic it follows directly from Theorem \ref{main} and Lemma
\ref{vanish}. Since the Borel-Serre compactification exists for arbitrary
Fuchsian groups, the proof runs through and we also get Theorem
\ref{main-Fuchs} for torsion-free Fuchsian groups.
An arbitrary Fuchsian group $\Ga$ has a finite index subgroup $\Ga'$ which
is torsion-free. An inspection shows that all our constructions allow
descent from $\Ga'$-invariants to $\Ga$-invariants and thus Theorem
\ref{main-Fuchs} follows.
\qed

%%--------------------Here the manuscript ends--------------------------------
\Addresses

\newpage
\begin{thebibliography}{XXX}

\bibitem{Arthur}
 \bf Arthur, J.:
 \it Eisenstein series and the trace formula.
 \rm Automorphic Forms, Representations, and L-Functions. Corvallis 1977; Proc. Symp. Pure Math. XXXIII, 253-274
(1979).

\bibitem{BoGa}
\bf Borel, A.; Garland, H.:
 \it Laplacian and discrete spectrum of an arithmetic group.
 \rm Amer.\ J.\ Math. 105, 309--335 (1983). 

\bibitem{BS}
 \bf Borel, A.; Serre, J.P.:
 \it Corners and Arithmetic Groups.
 \rm Comment.\ Math.\ Helv. 48, 436--491 (1973).

\bibitem{Borel-stableI}
 \bf Borel, A.:
 \it Stable real cohomology of arithmetic groups.
 \rm Ann. Sci. ENS 7, 235-272 (1974).

\bibitem{BorWall}
 \bf Borel, A.; Wallach,N.:
 \it Continuous Cohomology, Discrete Groups, and Representations of Reductive Groups.
 \rm Ann. Math. Stud. 94, Princeton 1980.

\bibitem{Borel-stableII}
\bf Borel, A.:
\it Stable real cohomology of arithmetic groups. II.
\rm Manifolds and Lie groups (Notre Dame, Ind., 1980), pp. 21--55, Progr. Math.,
14, Birkh\"auser, Boston, Mass., 1981.

\bibitem{BO1}
 \bf Bunke, U.; Olbrich, M.:
 \it $\Gamma$-Cohomology and the Selberg Zeta Function
 \rm J. reine u. angew. Math. 467, 199-219 (1995).

\bibitem{BO2}
\bf Bunke, U.; Olbrich, M.:
\it  Cohomological properties of the canonical globalizations of Harish-Chandra
modules. Consequences of theorems of Kashiwara-Schmid, Casselman, and
Schneider-Stuhler.
\rm Ann. Global Anal. Geom. 15
(1997), no. 5, 401--418.

\bibitem{BO3}
\bf Bunke, U.; Olbrich, M.:
\it Resolutions of distribution globalizations of Harish-Chandra modules and cohomology. 
\rm J. Reine Angew. Math. 497, 47-81, (1998).

\bibitem{C}
 \bf Casselman, W.:
 \it Canonical extensions of Harish-Chandra modules to representations of G.
 \rm Can. J. Math. 41, 385-438 (1989).

\bibitem{D1}
 \bf Deitmar, A.:
 \it Geometric zeta-functions of locally symmetric spaces.
 \rm Am. J. Math. 122, vol 5, 887-926 (2000).

\bibitem{D2}
\bf Deitmar, A.:
\it Selberg zeta functions for spaces of higher rank.
\rm http://arXiv.org/abs/math.NT/0209383.

\bibitem{F}
\bf Franke, J.:
\it Harmonic analysis in weighted $L\sb 2$-spaces. 
\rm  Ann. Sci. \'Ecole Norm. Sup. (4) 31 (1998), no. 2, 181--279.

\bibitem{Gro}
 \bf Grothendieck, A.:
 \it Produits tensoriels topologiques et espaces nucl\'eaires.
 \rm Mem. Amer. Math. Soc.  1955  (1955), no. 16.

\bibitem{HeSch}
 \bf Hecht, H.; Schmid, W.:
 \it Characters, asymptotics and $\n$-homology of Harish-Chandra modules.
 \rm Acta Math. 151, 49-151 (1983).

\bibitem{I}
\bf Iwaniec, H.:
\it Spectral methods of automorphic forms. Second edition.
\rm Graduate Studies in Mathematics, 53. 
American Mathematical Society, Providence, RI; Revista Matem\'atica
Iberoamericana, Madrid, 2002.

\bibitem{Juhl}
 \bf Juhl, A.:
 \it Cohomological theory of dynamical zeta functions.
 \rm Progress in Mathematics, 194. Birkhäuser Verlag, Basel,
2001.

\bibitem{KS}
\name{Kashiwara, M.; Schmid, W.:} \it Quasi-equivariant
{\cal D}-modules, equivariant derived category, and
representations of reductive Lie groups. \rm In: Lie
Theory and Geometry. In Honor of Bertram Kostant, Progr.
in Math., Birkh\"auser, Boston. 457-488 (1994).

\bibitem{L}
\bf Lang, S.:
\it Algebra. 
\rm Revised third edition. Graduate Texts in
Mathematics, 211. Springer-Verlag, New York, 2002.

\bibitem{Langlands}
 \bf Langlands, R.:
 \it On the Functional Equations Satisfied by Eisenstein Series.
 \rm SLNM 544, 1976.

\bibitem{LZ}
\bf Lewis, J.; Zagier, D.:
\it Period functions for Maass wave forms.
\rm I. Ann. of
Math. (2) 153 (2001), no. 1, 191--258.

\bibitem{Sch}
\bf Schaefer, H.:
\it Topological vector spaces.
\rm The Macmillan Co.,
New York; Collier-Macmillan Ltd., London 1966

\bibitem{S}
\bf Schmidt, T.:
\it Rational Period Functions and Parabolic Cohomology.
\rm J. Number Theory 57, 50-65 (1996).

\bibitem{Tr}
\bf Treves, F.:
\it Topological vector spaces, distributions and kenels.
\rm Academic Press, New York 1967.


\bibitem{Zagier}
\bf Zagier, D.:
\it New Points of View on the Selberg Zeta Function.
\rm Proceedings of the Japanese-German Seminar ``Explicit structures of modular
forms and zeta functions'' Hakuba, Sept. 2001. Ryushi-do, 2002.

\bibitem{Zimmer}
 \bf Zimmer, R.:
 \it Ergodic Theory and Semisimple Groups
 \rm Birkh\"auser 1984.

\end{thebibliography}
\end{document}